\documentclass{article}

\usepackage{amsmath}
\usepackage{amssymb}

\newtheorem{thm}{Theorem}

\title{Shortest Distance in Modular Cubic Polynomials}
\author{Tsz Ho Chan}
\date{}

\begin{document}
\maketitle

\begin{abstract}
In this paper, we study how small a box contains at least two points from a modular cubic polynomial $y \equiv a x^3 + b x^2 + c x + d \pmod p$ with $(a, p) = 1$. We prove that some square of side length $p^{1/6 + \epsilon}$ contains two such points.
\end{abstract}

\section{Introduction and Main results}
Recently, the author \cite{C} studied the shortest distance in modular hyperbola $x y \equiv c \pmod p$ and its relation with the least quadratic nonresidue modulo $p$. Inspired by this, we study the shortest distance in modular cubic polynomial in this paper. It is worthwhile to remark that, the shortest distance can have magnitude $p$ for linear polynomials while the shortest distance can be as small as $O(1)$ for quadratic polynomials.

\bigskip

Let $p > 3$ be a prime, $(a, p) = 1$ and any integer $c$. We consider the modular reduced cubic
\[
C_{a, c} := \{ (x, y): y \equiv a x^3 + c x \pmod p \}.
\]
The restriction to such reduced cubic polynomials is not restrictive at all as one can transform a general cubic to such form through change of variables in $x$ and $y$ which does not affect the shortest distance between two points in $C_{a, c}$. Instead of distances, we consider how small a box
\[
B(X, Y; H) := \{ (x, y): X + 1 \le x \le X + H \pmod p, Y + 1 \le y \le Y + H \pmod p \}
\]
contains at least two points in $C_{a, c}$ where $X$ and $Y$ run over $0$, $1$, ..., $p-1$.

To study this, we use a recent result of Heath-Brown \cite{H} and Shao \cite{S} on mean-value estimates of character sums:
\begin{thm} \label{Shao}
Given $H \le p$, a positive integer and any $\epsilon > 0$. Suppose that $0 \le N_1 < N_2 < ... < N_J < p$ are integers satisfying $N_{j+1} - N_j \ge H$ for $1 \le j < J$. Then
\[
\sum_{j = 1}^{J} \max_{h \le H} |S(N_j; h)|^{2r} \ll_{\epsilon, r} H^{2r-2} p^{1/2 + 1/(2r) + \epsilon}
\]
where
\[
S(N; H) := \sum_{N < n \le N + H} \chi(n)
\]
and $\chi$ is any non-principal character modulo $p$.
\end{thm}

Applying the above theorem, we can show that
\begin{thm} \label{thm1}
For any $\epsilon  > 0$, for any $(a, p) = 1$, integer $c$ and $H \gg_\epsilon p^{1/6 + \epsilon}$, we have
\[
| C_{a,c} \cap B(X, Y; H) | \ge 2
\]
for some $0 \le X , Y \le p-1$.
\end{thm}

\bigskip

{\bf Some Notations} Throughout the paper, $p$ stands for a prime. The symbol $|S|$ denotes the number of elements in the set $S$. We also use the Legendre symbol $(\frac{\cdot}{p})$. The notations $f(x) \ll g(x)$, $g(x) \gg f(x)$ and $f(x) = O(g(x))$ are equivalent to $|f(x)| \leq C g(x)$ for some constant $C > 0$. Finally, $f(x) \ll_{\lambda_1, ..., \lambda_k} g(x)$, $g(x) \gg_{\lambda_1, ..., \lambda_k} f(x)$ and $f(x) = O_{\lambda_1, ..., \lambda_k} (g(x))$ mean that the implicit constant $C$ may depend on $\lambda_1$, ..., $\lambda_k$.

%-------------------------------------------------------------------------------------------
\section{The Basic Argument} \label{basic}
Without loss of generality, we assume that $p > 3$. For $(a, p) = 1$ and any integer $c$, suppose $| C_{a,c} \cap B(X, Y; H) | \ge 2$ for some $0 \le X, Y \le p-1$. This means that
\begin{equation} \label{start}
y \equiv a x^3 + c x \pmod p, \text{ and } y + v \equiv a (x+u)^3 + c(x+u) \pmod p
\end{equation}
for some $1 \le x, y \le p$ and $1 \le u, v \le H$. Subtracting, we get
\[
v \equiv 3 a u (x^2 + u x + \overline{3}u^2) + c u \pmod p
\]
where $\overline{y}$ denotes the multiplicative inverse of $y$ modulo $p$ (i.e. $y \overline{y} \equiv 1 \pmod p$.)
After some algebra and completing the square, we have
\[
(2x + u)^2 \equiv 4 \overline{3} v \overline{a} \overline{u} - \overline{3} u^2 - 4 \overline{3} \overline{a} c \pmod p.
\]
The above process is reversible. So $| C_{a,c} \cap B(X, Y; H) | \ge 2$ for some $0 \le X, Y \le p-1$ is equivalent to
\[
\Bigl( \frac{-3}{p} \Bigr) \Bigl( \frac{a}{p} \Bigr) \Bigl( \frac{u}{p} \Bigr) \Bigl( \frac{a u^3 + 4 c u - 4 v}{p} \Bigr) = 1.
\]
We are going to restrict our attention to even $u = 2u'$'s and $v = 2v'$'s. So we want
\begin{equation} \label{cond}
\Bigl( \frac{-3}{p} \Bigr) \Bigl( \frac{a}{p} \Bigr) \Bigl( \frac{u'}{p} \Bigr) \Bigl( \frac{a u'^3 + c u' - v'}{p} \Bigr) = 1 \text{ for some } 1 \le u', v' \le H/2.
\end{equation}

%-------------------------------------------------------------------------------------------
\section{Proof of Theorem \ref{thm1}}
Suppose (\ref{cond}) is not true. Then either
\[
\Bigl( \frac{-3}{p} \Bigr) \Bigl( \frac{a}{p} \Bigr) \Bigl( \frac{u'}{p} \Bigr) \Bigl( \frac{a u'^3 + c u' - v'}{p} \Bigr) = 0;
\]
or
\[
\Bigl( \frac{-3}{p} \Bigr) \Bigl( \frac{a}{p} \Bigr) \Bigl( \frac{u'}{p} \Bigr) \Bigl( \frac{a u'^3 + c u' - v'}{p} \Bigr) = -1
\]
for all $1 \le u', v' \le H/2$. If the former is true for two pairs of $1 \le u', v' \le H/2$, we have
\begin{equation} \label{aside}
a u'^3 + c u' \equiv v' \pmod p \text{ and } a u''^3 + c u'' \equiv v'' \pmod p
\end{equation}
which gives Theorem \ref{thm1}. Henceforth we suppose the latter is true for all but at most one pair of $1 \le u', v' \le H/2$. Hence
\begin{align*}
H^2 \ll& \Big| \sum_{u' \le H/2} \sum_{v' \le H/2} \Bigl( \frac{u'}{p} \Bigr) \Bigl( \frac{a u'^3 + c u' - v'}{p} \Bigr) \Big| \le \sum_{u' \le H/2} \Big| \sum_{v' \le H/2} \Bigl( \frac{a u'^3 + c u' - v'}{p} \Bigr) \Big| \\
\le& \Bigl( \sum_{u' \le H/2} 1 \Bigr)^{(2r - 1)/(2r)} \Bigl( \sum_{u' \le H/2} \Big| \sum_{v' \le H/2} \Bigl( \frac{a u'^3 + c u' - v'}{p} \Bigr) \Big|^{2r} \Bigr)^{1/(2r)}.
\end{align*}
Suppose $| C_{a,c} \cap B(X, Y; H) | \le 1$ for all $0 \le X, Y \le p-1$. Then the points $a u'^3 + c u'$ are spaced more than $H$ apart. So we can apply Theorem \ref{Shao} and get
\[
H^2 \ll_{\epsilon, r} H^{(2r - 1)/(2r)} (H^{2r - 2} p^{1/2 + 1/(2r) + \epsilon})^{1/(2r)}
\]
which gives $H \ll_{\epsilon, r} p^{(r+1)/(6r) + \epsilon/2}$. This contradicts $H \gg_\epsilon p^{1/6 + \epsilon}$ if $r$ is sufficiently large. This final contradiction together with (\ref{aside}) gives Theorem \ref{thm1}.

%-------------------------------------------------------------------------------------------

Tsz Ho Chan \\
Department of Mathematical Sciences \\
University of Memphis \\
Memphis, TN 38152 \\
U.S.A. \\
thchan6174@gmail.com \\

\end{document}